\documentclass{amsart}
\usepackage{latexsym,amssymb,amsmath,amsthm}
\usepackage[english]{babel}

\newcommand{\atc}{\`a}

\newcommand{\R}{{\mathbb{R}}}
\newcommand{\C}{{\mathbb{C}}}

\newcommand{\n}{\noindent}
\newcommand{\mioplus}[1][]{\,\,\oplus_{#1}\,\,}
\def\P{\text{$\mathbb{P}$}}
\def\J{\text{$\mathbb{J}$}}
\def\H{\text{$\mathbb{H}$}}
\def\Z{\text{$\mathbb{Z}$}}
\newcommand\SO{{\rm SO}}

\newcommand\SU{{\rm SU}}
\def\Sp{{\rm Sp}}

\newcommand\U{{\rm U}}
\newcommand\SUr[2]{{{\rm S}({\rm U}_{#1}\times {\rm U}_{#2})}}
\newcommand\GL{{\rm G \ell}}

\newcommand\Irr{{\rm Irr}}
\newcommand\I{{\rm I}}
\newcommand\T{{\rm T}}

\newcommand\mf[1]{\mathfrak{#1}}
\newcommand\Spin{{\rm Spin }}

\newcommand\Ea{{\rm E}}
\newcommand\Fa{{\rm F}}

\newcommand\K{{\rm K}}
\newcommand\Se{{\rm S}}
\newcommand\li[1]{\mathfrak{l}_{#1}}

\newcommand\su[1]{\mathfrak{su}(#1) }
\newcommand\so[1]{\mathfrak{so} (#1) }
\newcommand\un[1]{\mathfrak{u} (#1) }
\newcommand\lra{\longrightarrow}
\newcommand\asp{\mathfrak{sp}}

\newcommand\Aut{{\rm Aut}}

\newtheorem{thm}{Theorem}[section]
\newtheorem{prop}{Proposition}[section]
\newtheorem{lem}{Lemma}[section]

\newtheorem{con}{Conjecture}

\theoremstyle{definition}

\theoremstyle{remark}

\numberwithin{equation}{section}

\begin{document}
\title[Coisotropic and Polar actions]{Coisotropic and polar actions on
compact irreducible Hermitian symmetric spaces}
\author{Leonardo Biliotti}
\curraddr{Dipartimento di Matematica, Universit\atc $\ $
Politecnica delle Marche, Via Brecce Bianche, 60131, Ancona,
Italy} \email{biliotti@dipmat.univpm.it}
\begin{abstract}
We obtain the full classification of coisotropic and polar
isometric actions of compact Lie groups on irreducible Hermitian
symmetric spaces.
\end{abstract}
\maketitle
%
%
\section{Introduction}
The aim of the present paper is to investigate polar and
coisotropic actions on compact irreducible Hermitian symmetric
spaces.

The action of a compact Lie group $K$ of isometries on a
Riemannian manifold $(M,g)$ is called {\em polar} if there exists
a properly embedded submanifold $\Sigma$ which meets every
$K-$orbit and is orthogonal to the $K-$orbits in all common
points. Such a submanifold $\Sigma$ is called a {\em section} (see
\cite{PT1}, \cite{PT2}) and it is automatically totally geodesic;
if it is flat, the action is called {\em hyperpolar}.

Let $(M,g)$ be a compact K\"ahler manifold with K\"ahler form
$\omega$ and let $K$ be a compact connected Lie subgroup of its
full isometry group. The $K$-action is called {\em coisotropic} or
{\em multiplicity free} if the principal $K$-orbits are
coisotropic with respect to $\omega$ \cite{Hw}. Notice that the
existence of an open subset consisting of coisotropic orbits
implies that all $K-$orbits are coisotropic, see \cite{Hw}.
Multiplicity free representations form a very restricted class of
representation. Nevertheless they are very important since every
``nice'' result in the invariant theory of particular
representations can be traced back to a multiplicity free
representation. This holds for example for a Capelli identities
\cite{HM} and also all of Weyl's first and second fundamental
theorems can  be explained by some multiplicity freeness result.

Kac \cite{Kac} and Benson and Ratcliff \cite{BR} have given  the
classification of linear multiplicity free representations, from
which one has the full classification of coisotropic actions on
$Gr (k,n)$ for $k=1$, i.e. on the complex projective space. In a
recent paper \cite{BA} the complete classification of polar and
coisotropic actions on complex Grassmannians has been obtained
while in \cite{PoT}, as an application of the main result, it was
given the complete classification of this kind of actions on the
quadric $\SO (n+2) / \SO(2) \times \SO(n)$. Hence it is natural to
investigate coisotropic and polar actions on the other compact
irreducible Hermitian symmetric spaces, which are $\SO(2m) /
\U(m),$ $\Sp(m) / \U (m),$ $\Ea_7 / \T^1 \cdot \Ea_6$ and $\Ea_6 /
\T^1 \cdot \Spin (10).$ Our main result is given in the following
%
%
\begin{thm} \label{risultato} Let $K$ be a compact
connected Lie subgroup of $\Sp(m),$ respectively $\SO (2m),$
acting non-transitively on the Hermitian symmetric space $M=
\Sp(2m)/ \U (m),$ respectively $M=\SO (2m) / \U(m).$ Then $K$ acts
coisotropically on $M$ if and only if its Lie algebra $\mf{k},$ up
to conjugation in $\asp (2m),$ respectively $\mf{o} (2m),$
contains one of the Lie algebras  appearing in Table $1.$ In Table
$2$ we list, up to conjugation, all the subgroups of $\Ea_7,\
\Ea_6,$ which act non-transitively and coisotropically on $\Ea_7 /
\T^1 \cdot \Ea_6 $ and $\Ea_6 / \T^1 \cdot \Spin (10) $
respectively.
\end{thm}
\begin{small}
\begin{table}[ht]
\caption{}\label{eqtable}
\renewcommand\arraystretch{1.0}
\noindent\[
\begin{array}{|c|c|c|}\hline
\mf k & M & {\rm conditions}\\ \hline\hline \un 1 & \Sp (1) /
\U(1) & \\ \hline \su m & \Sp (m) / \U(m) & m \geq 2 \\ \hline
\asp(k) + \asp(m-k) & \Sp(m) / \U(m) & 1 \leq k \leq m-1\\ \hline
\asp(m-1) + \un 1 & \Sp(m) / \U(m) & m \geq 2 \\ \hline \asp (m) +
\asp(1) +\asp(1) & \Sp(m+2) / \U(m+2) & \\ \hline\hline
%
%
%
%
%
\R (0)      & \SO (4) / \U(2) & \R(0) \ \mathrm{line\ in} \ \mf
t_2 \times \mf z
                                     \\ \hline
\mf z + \mf t_3 & \SO (6) / \U(3) &  \\ \hline
\R (\frac{1}{2k} ) + \su{2k}   & \SO (4k+2) / \U (2k+1)  &  k\geq
2,\ \R(\frac{1}{2k}) \ \mathrm{line \ in} \ \mf a \times \mf z \\
\hline \R  + \su{2k+1}   & \SO (4k+4) / \U (2k+2)  & k\geq 2,\ \R
\ \mathrm{means \ any \ line \ in} \ \mf a \times \mf z \\ \hline
\R (0) + \su 3 & \SO (8) / \U (4) & \R(0) \ \mathrm{line \ in} \
\mf a \times \mf z
\\ \hline
\mf z + \su 2   & \SO (6) / \U (3) & \\ \hline
\su m       & \SO(m) / \U(m) & m \geq 2 \\ \hline
\mf z + \asp (2)  & \SO (8) / \U (4)  &  \\ \hline \asp (1) +
\asp(2) & \SO(8) / \U(4) & \asp(1) \otimes \asp(2) \subseteq \so 8
\\ \hline
%
%
%
\mf{so} (k) + \mf{so}(2m-k) & \SO(2m) / \U (m) & \\ \hline
\mf{so}(2m-2) & \SO(2m) / \U(m) & m \geq 3 \\ \hline
\mf{so}(2m-6) + \mf{u} (3) & \SO(2m)/ \U(m) & m \geq 5  \\ \hline
\mf{so}(2m-4) + \un 2 & \SO(2m)/ \U(m) & m \geq 4\\ \hline
\mf{so}(2m) + \R(1,-1) & \SO(2(m+2)) / \U(m+2) & m \geq 5,\\
& & \R(1,-1)\subseteq \so 2 \times \so 2 \subseteq \so 4
\\ \hline
\so{4} + \so 2 + \so 2 & \SO(8) / \U(4) & \\ \hline
\mf{g}_2 & \SO(8)/ \U(4) & \mf{g}_2 \subset \so 7 \subset \so 8 \\
\hline
\end{array}
\]
\end{table}
\end{small}
%
%
\begin{small}
\begin{table}[ht]
\caption{}\label{eqtable}
\renewcommand\arraystretch{1.0}
\noindent\[
\begin{array}{|c|c|c|c|c|}\hline
&  M= \Ea_7 / \T^1 \cdot \Ea_6  & & &  \\ \hline {\rm maximal \
subgroups} & \T^1 \cdot \Ea_6  & \SU(2) \cdot \Spin (12) & \SU(8)
/ \Z_2  &  \\ \hline & & \T^1 \cdot \Spin (12)    & \SUr{1}{7}/
\Z_2 & \\ \hline & & \SU(2) \cdot \Spin(11) &   \SU(7) / \Z_2
&  \\
\hline\hline
%
%
&   M= \Ea_6 / \T^1 \cdot \Spin (10)  & & &  \\ \hline {\rm
maximal \ subgroups} & \T^1 \cdot \Spin (10) & \Sp(1) \cdot \SU(6)
& \Sp(4)/ \Z_2                 & {\rm F}_4                 \\
\hline & \Spin (10)             &  \T^1 \cdot \SU(6)    &  & \\
\hline & \T^1 \cdot \Spin (9)  &  \Sp (1) \cdot \U(5)   &  &
\\
\hline & \T^1 \cdot (\T^1 \times \Spin (8)) & \T^1 \cdot \U(5) & &
\\ \hline
\end{array}
\]
\end{table}
\end{small}

All the Lie algebras listed in the first column, unless explicitly
specified, are meant to be standardly embedded into $\asp (m),$
respectively $\mf{so}(2m)$, e.g. $\asp (m) + \mf{u} (1)  \subset
\asp (m)+ \asp (2) \subset \asp (m),$ $\mf{so}(2m-3)+ \mf u (3)
\subset \mf{so} (2m-3) + \mf{so}(6) \subset \mf{so}(2m).$ The
notations used in Table $1$ are as follows. We denote with
$\mf{z}$ the one dimensional center of $Lie(\U(m)),$ with $\mf a$
the centralizer of the semisimple part of $\mf k$ in $\su m
\subset$ $Lie( \U(m))$ and with $\mf t_m$ the subalgebra of a
maximal torus of $\SU (m) \subset \U(m)$ With this notation
$\R(\alpha)$ denotes any line in $\mf a \times \mf z$ different
from $y=\alpha x$ while $\R(1,-1) \subseteq \so 2 + \so 2
\subseteq \so 4$ means any line in the plane $\so 2 \times \so 2$
different from $y=x$ and $y=-x.$ Finally, in Table $2$ the
juxtaposition $A \cdot B$ of two groups generally denotes the
quotient $A \times_{\Z_2} B.$

Victor Kac \cite{Kac} obtained a complete classification (Tables
Ia, Ib, in the Appendix) of irreducible multiplicity free actions
$(\sigma,V)$. Most of these include a copy of the scalars $\C$
acting on $V$. We will say that a multiplicity free action
$(\sigma,V)$ of a complex group $G$ is {\em decomposable} if we
can write $V$ as the direct sum $V=V_1\oplus V_2$ of proper
$\sigma(G)$-invariant subspaces in such a way that
$\sigma(G)=\sigma_1(G)\times\sigma_2(G)$, where $\sigma_i$ denotes
the restriction of $\sigma$ to $V_i$. If $V$ does not admit such a
decomposition then we say that $(\sigma,V)$ is an {\em
indecomposable} multiplicity free action. C. Benson and G.
Ratcliff have given the complete classification of indecomposable
multiplicity free actions (Tables IIa, IIb in the Appendix). We
recall here their theorem (Theorem 2, page 154 \cite{BR})
\begin{thm}\label{Ben}
Let $(\sigma,V)$ be a regular representation of a connected
semisimple  complex algebraic group $G$ and decompose $V$ as a
direct sum of $\sigma(G)$-irreducible subspaces, $V=V_1\oplus V_2
\oplus \cdots \oplus V_r$. The action of $(\C^*)^r \times G$ on
$V$ is an indecomposable multiplicity free action if and only if
either
\begin{description}
\item[(1)] $r=1$ and $\sigma(G)\subseteq \GL(V)$ appears in Table
Ia (see the Appendix); \item[(2)] $r=2$ and $\sigma(G)\subseteq
\GL (V_1)\times \GL( V_2)$ appears in Tables IIa and IIb (see the
Appendix).
\end{description}
\end{thm}
In \cite{BR} are also given conditions under which one can {\em
remove} or {\em reduce} the copies of the scalars preserving the
multiplicity free action. Obviously  if an action is coisotropic
it continues to be coisotropic also when this action includes
another copy of the scalars. We will call {\em minimal} those
coisotropic actions in which the scalars, if they appear, cannot
be reduced.

Let $K$ be a compact group acting isometrically on a compact
K\"ahler manifold $M.$ This action is automatically holomorphic by
a theorem of Kostant ( see \cite{KN}, vol I, page 247) and it
induces by compactness of $M$ an action of the complexified group
$K^{\C}$ on $M.$ We say that $M$ is $K^\C$-{\em almost
homogeneous} if $K^\C$ has an open orbit in $M.$ If all Borel
subgroups of $K^\C$ act with an open orbit on $M$, then the
$K^\C$-open orbit $\Omega$ is called a {\em spherical homogeneous
space} and $M$ is called a {\em spherical embedding} of $\Omega.$
We will briefly recall some results that will be used in the
sequel.
\begin{thm}{\label{class}}\cite{Hw}
Let $M$ be a connected compact K\"ahler manifold with an isometric
action of a connected compact group $K$ that is also Poisson. Then
the following conditions are equivalent:
\begin{itemize}
\item [(i)] The $K$-action is coisotropic. \item [(ii)] The
cohomogeneity of the $K$ action is equal to the difference between
the rank of $K$ and the rank of a regular isotropy subgroup of
$K$. \item [(iii)] The moment map
$\mu\,:\,M\rightarrow{\mathfrak{k}^*}$ separates orbits. \item
[(iv)] The K\"ahler manifold $M$ is projective algebraic,
$K^\mathbb{C}$-almost homogeneous and a spherical embedding of the
open $K^\mathbb{C}$-orbit.
\end{itemize}
\end{thm}
We remark here that conditions (i) to (iii) are equivalent even
without the hypothesis of compactness on $M$ (see \cite{Hw}).

As an immediate consequence of the above theorem one can deduce,
under the same hypotheses on $K$ and $M,$ two simple facts that
will be frequently used in our classification:
\begin{itemize}
\item [1] Let $p$ be a fixed point on $M$ for the $K$-action, or
$Kp$ a complex $K$-orbit, then the $K$-action is coisotropic if
and only if the slice representation is coisotropic (see \cite{Hw}
page 274). \item [2] {\em dimensional condition.} If $K$ acts
coisotropically on $M$ the dimension of a Borel subgroup $B$ of
$K^\C$ is not less than the dimension of $M$.
\end{itemize}
A relatively large class of coisotropic actions is provided by
polar ones. A result due to Hermann (\cite{He}) states that given
$K$  a compact Lie group and two symmetric subgroups
$H_1$,$H_2\subseteq{K}$,  then $H_i$ acts hyperpolarly on $K / H_j
$ for $i,{j}\in{1,2}.$ This kinds of action are coisotropic since
for \cite{PoT} a  polar action on an irreducible compact
homogeneous K\"ahler manifold is co\-iso\-tro\-pic.
%
%

Once we shall determined the complete list of coisotropic actions
on compact irreducible Hermitian symmetric spaces we have also
investigated which ones are polar. Dadok \cite{Da}, Heintze and
Eschenburg \cite{He} have classified the irreducible polar linear
representations, while I.Bergmann \cite{Bergmann} has found all
the reducible ones. Using their results we determine in section
$7$ the complete classification of the polar actions on the
following Hermitian symmetric spaces $\SO(2m)/ \U(m),$ $\Sp(m) /
\U(m),$ $\Ea_6 / \T^1 \cdot \Spin (10),$ $\Ea_7 / \T^1 \cdot \Ea_6
.$ An interesting consequence of this classification is that the
polar actions on these manifolds  are just the hyperpolar ones.
The same result holds  on the quadrics (see \cite{PoT}) and on the
complex Grassmannians (see \cite{BA}). In particular, we have the
following
\begin{prop} \label{pi}
A polar action on a compact irreducible Hermitian symmetric space
of rank bigger than one is hyperpolar.
\end{prop}
This is in contrast to complex projective space or more generally
to rank one symmetric spaces that admit many polar actions that
are not hyperpolar (see \cite{PoT2}).

We point out also that on the Hermitian symmetric space $M=\Ea_7 /
\T^1 \cdot \Ea_6,$ respectively $M= \Sp(m) / \U(m),$ our result
implies that a compact connected Lie subgroup $K$ of $\Ea_7,$
respectively $\Sp(m),$ acts polarly on $M$ if and only if $K$ is a
symmetric group.

We mention the following conjecture concerning the nature of polar
actions on compact symmetric spaces.
\begin{con}
A polar action on a  compact symmetric space of rank bigger than
one is hyperpolar.
\end{con}
In particular in Proposition \ref{pi} is given  the positive
answer in the class of compact irreducible Hermitian symmetric
spaces.

The complete classification of polar actions on the compact
irreducible Hermitian symmetric spaces, which they have been
investigated in this paper, is given in the following
%
%
%
\begin{thm} \label{polar}
Let $K$ be a compact connected Lie subgroup of $\SO(2m),$
respectively $\Sp(m),$ acting non-transitively on $M=\SO(2m)/ \U
(m)$ respectively $M=\Sp(m) / \U(m) $. Then $K$ acts polarly on
$M$ if and only if its Lie algebra $\mf{k}$ is  conjugate, in $\mf
o (2m),$ respectively $\asp (m),$ to one of the Lie algebras
appearing in Table $3.$ In Table $4$ we list, up to conjugation,
all the subgroups of $\Ea_7,\ \Ea_6 ,$ which act non-transitively
and polarly on $\Ea_7/ \T^1 \cdot \Ea_6$ and $\Ea_6 / \T^1 \cdot
\Spin (10)$ respectively. In particular on these manifolds is that
polar actions are hyperpolar.
\end{thm}
\newpage
\begin{small}
\begin{table}[ht]
\caption{}\label{eqtable}
\renewcommand\arraystretch{1.0}
\noindent\[
%
%
%
\begin{array}{|c|c|c|}\hline \mf k & M & {\rm conditions}\\
\hline\hline
\mf u (m) & \Sp (m) / \U(m) & m \geq 1 \\ \hline
\mf{sp}(k) + \mf{sp}(m-k) & \Sp(2m) / \U(m) & \\ \hline\hline
%
\un m & \SO(2m) / \U(m) & \\ \hline
\su m & \SO(2m) / \U(m) & m \  {\rm odd} \\ \hline
%
%
\mf{so}(k) + \mf{so}(2m-k) & \SO(2m) / \U(m) &   \\  \hline
\mf{so}(2m-2) & \SO(2m)/ \U(m) & m \geq 3 \\ \hline
\mf{g}_2  & \SO(8) / \U(8) & \mf{g}_2 \subset \mf{so}(7) \subset
\mf{so}(8) \\ \hline
\R(0) & \SO(4) / \U(2) & \\ \hline
\end{array}
\]
\end{table}
\end{small}
%
%
%
%
\begin{small}
\begin{table}[ht]
\caption{}\label{eqtable}
\renewcommand\arraystretch{1.0}
\noindent\[
%
%
%
\begin{array}{|c|c|c|c|}\hline & M=\Ea_7 / \T^1 \cdot \Ea_6 & &
\\ \hline
\T^1 \cdot \Ea_6 & \Spin (12) \cdot \SU (2) & \SU(8) / \Z_2 & \\
\hline\hline
& M=\Ea_6/ \T^1 \cdot \Spin (10)  &  & \\ \hline
\T^1 \cdot \Spin(10)  & \SU(8) / \Z_2 & \Sp(4)/ \Z_2  & {\rm F}_4
\\ \hline
\Spin (10)            &               &            &   \\ \hline
\end{array}
\]
\end{table}
\end{small}
%
%
%
%
We here briefly explain our method in order to prove our main
theorem. Thanks to Theorem \ref{class}, (iv) we have that if $K$
is a subgroup of a compact Lie group $L$ such that $K$ acts
coisotropically on $M$ so does $L$. As a consequence, in order to
classify coisotropic actions on  $\SO(2m)/ \U (m)$ ($\Sp(m) /
\U(m),$ $\Ea_7 / \T^1 \cdot \Ea_6,$ $\Ea_6 / \T^1 \cdot \Spin
(10)$), one  may suggest a sort of ``telescopic'' procedure by
restricting to maximal subgroups $K$ of $\SO(n)$,($\Sp(m),$
$\Ea_7$ $\Ea_6$) hence passing to maximal subgroups that give rise
to coisotropic actions and so on.

This paper is organized as follows. In section 2 we prove a useful
result that we shall use throughout this paper. From section 3 to
section 6 we give the proof of Theorem \ref{risultato}. We have
divided every section in subsections in each of which we analyze
separately one of the maximal subgroups of $\SO(m)$ respectively
$\Sp (2m),$ $\Ea_7$ and $\Ea_6.$ In the seventh section we give
the proof of Theorem \ref{polar}.

We enclose, in the Appendix, the tables of irreducible  and
reducible linear multiplicity free representations (Tables Ia, Ib
and Tables IIa, IIb respectively) and the tables of maximal
subgroups of $\Sp(2m),$ $ \SO(n)$ and $\SU(n)$ (Tables {\bf III},
{\bf IV},{\bf V}).
%
%
\section{Preliminaries}
Let $\mf g$ be a Lie semisimple complex algebra. We will denote by
$\mf b$ a Borel Lie algebra of $\mf g,$ whose dimension is
$\frac{1}{2} (\dim \mf g + r (\mf g ) ),$ where $r (\mf g)$ is the
dimension of a Cartan subalgebra, namely the rank of $\mf g.$
Throughout this paper we will identify the fundamental dominant
weights $\Lambda_l$ with the corresponding irreducible
representations. It is well known that any irreducible
representation corresponds to a highest weight $\sigma$ and any
highest weight is of the form $\sigma= \sum_i m_i \Lambda_i,$
where $m_i$ are non-negative integers. We will denote by
$d(\sigma)$ the representation degree of $\sigma,$ i.e. the
dimension of the vector space on which $\mf g$ acts with the
irreducible representation $\sigma.$ Using the Weyl's dimensional
formula it easy to check that if $m_i \geq n_i$ then $d (\sum_i
m_i \Lambda_i) \geq d (\sum_i n_i \Lambda_i)$ and the equality
hold if and only if $m_i=n_i.$
\begin{lem} \label{zio}
Let $\mf g$ be a simple complex Lie algebra and let $\sigma: \mf g
\lra \mf{gl} (V)$ be a representation of $\mf g$ on $V$ with
$d=\dim V.$ Let $\mf{b}$ be the  Lie algebra of a Borel subgroup
of $\mf g.$ Then we have
\begin{enumerate}
\item $1+ \dim \mf b < \frac{1}{2}d(d-1)$ except when $\mf g=
\mf{sl} (m)$ and either $\sigma=\Lambda_1$ or $\sigma=
\Lambda_{m-1},$ $\mf g= \mf{sl}(2)$ and $\sigma = 2 \Lambda_1,$
$\mf g =\mf{so} (5)$ and $\sigma= \Lambda_2$ ({\it
spin}-representation) and $\mf g = \mf{so} (6)$ and either
$\sigma=\Lambda_3$ or $\sigma =\Lambda_2$ ({\it
spin}-representations); \item  $1+ \dim \mf b < \frac{1}{2}
d(d+1)$ except when $\mf g = \mf{sl} (m)$ and either $\sigma =
\Lambda_1$ or $\sigma = \Lambda_{m-1};$
\end{enumerate}
\end{lem}
\begin{proof}
Since the second affirmation can be deduced easily from the first,
we shall prove only our first statement. Our basic references are
\cite{SK} and \cite{Kol} Appendix B.
%

Assume $\mf g=\mf{sl}(m).$  Then the dimension of the Borel
subalgebra is $\dim \mf b = \frac{1}{2} (m-1)(m+2).$ The cases
$m=2,3$ are easy to check. If $m\geq 4,$ we have $d(\Lambda_1 +
\Lambda_{m-1}) \geq m+\frac{3}{2},$ $d(2 \Lambda_1) \geq
m+\frac{3}{2}$ and $d(\Lambda_2) \geq m+\frac{3}{2}.$ In
particular, for every representation $\sigma \neq \Lambda_1,
\Lambda_{m-1},$ one may verify that $1+ \dim \mf b <
\frac{1}{4}(2m+3)(2m+1) \leq \frac{1}{2} d(\sigma) (d(\sigma)
-1).$
%
%
Assume $\mf g = \asp (m), m \geq 3.$ Since $d(\sigma) >4m >
d(\Lambda_1)=2m,$ when $\sigma \neq \Lambda_1,$ we have $1+ \dim
\mf b= 1+ m^2+m < \frac{1}{2} d(\sigma) (d (\sigma) -1),$ since
$1+ m^2 + m <m(2m-1)$ for $m \geq 3.$
%
%
If $\mf g = \mf{so}(2m+1),$ we distinguish the case $m \geq 4$ and
$m=2,3.$ When $m\geq 4,$ since $d(\sigma) \geq 2m-1,$ we have $1+
\dim \mf b =m^2 +m < \frac{1}{2}d(\sigma) ( d(\sigma) -1).$ If
$m=3,$ since $d(\sigma) \geq 7,$ one may prove that
$\frac{1}{2}d(\sigma)(d(\sigma) -1) > 13$ is verified for  every
$\sigma,$ while in the case $m=2$ we have that $\sigma =
\Lambda_2$ does not satisfy the above inequality.
%
%
The case  $\mf g = \mf{so} (2m),$ can be resolved as before.
Indeed, if $m \geq 4$ then  it is to check that $d (\sigma ) \geq
2m-1$ for every $\sigma.$ In particular $1+\dim b=1+ m^2 <
(2m-1)(m-1) \leq \frac{1}{2}d(\sigma) (d(\sigma)-1).$ If $m=3,$
since $d(2 \Lambda_1)=d( \Lambda_1 + \Lambda_2 )=20,$ $d(\Lambda_1
+ \Lambda_3)=15,$ and $d(2 \Lambda_2 ) = d(2 \Lambda_3)= 10,$ one
may prove that $10 < \frac{1}{2} d(\sigma) (d(\sigma)-1)$ except
for  $\sigma=\Lambda_i,\ i=2,3.$
%
If $\mf g$ is of type $\mf g_2$ ( $\mf f_4 ,$ $\mf e_6,$ $\mf
e_7,$ $\mf e_8$) it is well known that  minimal representation
degree is $7$ (respectively $26,$ $27,$ $56,$ $248$) and the
dimension of a Borel subalgebra is $8$ (respectively $31,$ $42,$
$70,$ $127$), then for any representation $\sigma$ we have $1+
\dim \mf b <  \frac{1}{2} d(\sigma) (d(\sigma)-1).$
\end{proof}
%
%
%
\section{$M= \Sp (m) / \U (m)$}
%
%
\subsection{The case $K =\rho (H)$, $H$ simple such that $\rho$ $\in Irr_{\H}
(H)$, $d(\rho)=2m$}

We briefly explain our notation, that we will use throughout this
paper. Let $H$ be a simple group. By $Irr_{\R}(H)$, $Irr_{\C}(H)$,
$Irr_{\H}(H)$ we denote the irreducible representation of $H$ of
 real, complex and quaternionic type, see \cite{BTD}, Chapter II,
 $\S  6$ .

By Table III in the Appendix, if $K$ is the image of an
irreducible quaternionic representation $\rho$ of a simple Lie
group $H$, i.e. $K=\rho(H)$ where $\rho \in Irr_{\H} (H),$ then
$K$ is a maximal group. In this section we analyze this case.

Let $H$ be a simple group. It is well known that if $\mf h_o$ is a
simple real algebra whose complexification $\mf h$ is simple, its
irreducible representations are the restrictions of (uniquely
determined) irreducible representation of $\mf h.$ Our idea is
very simple: we impose the  {\it dimensional condition}. By lemma
\ref{zio} we have only to consider $(\mf{sl} (2), \Lambda_1),$
which corresponds to $\SU(2) \subseteq \U(2) \subseteq \SO(4).$
This case will be studied in next section, since $\SU(2)$ has a
fixed point.
%
%
\subsection{The fixed point case $K=\U(m)$} \label{fisso}

$\U(m)$ has a fixed point and the slice is given by $\Se^2
(\C^m).$ By Tables Ia and Ib, the action is multiplicity free and
the scalar can be removed when $m \geq 2.$ We will now go through
the maximal subgroups of $\U (m).$  Let $L\subset \U(m)$ be such
that Lie($L$)= $\mf z +\li 1,$ where $\li 1$ is a maximal
subalgebra of $\su m$ (see Table {\bf V} in the Appendix). By
lemma \ref{zio} the dimensional condition is not satisfies for
(i), (ii) and (v) of Table {\bf V}. The same holds for $\li 1= \su
p + \su q .$ Indeed, the dimension of a Borel subalgebra of $(\mf
z + \li 1)^{\C} $ is $1+\frac{1}{2} ( (p-1)(p+2)+(q-1)(q+2)).$ The
inequality $1+ \frac{1}{2}( (p-1)(p+2)+(q-1)(q+2) ) <
\frac{1}{2}(pq(pq+1))$ is always satisfies, so the action fails to
be multiplicity free. Indeed, let $f(x)=x^2(q^2-1)+x(q-1) -q^2
-q+2.$ Then $f'(x)=2x(q^2-1) +q-1>0,$ for $x \geq 3$ and
$f(3)=9(q^2-1)+3(q-1)-q^2-q+2>0,$ since $q \geq 2.$ Finally, if
$\li 1= \R + \su k + \su{m-k}$ then the slice becomes $\Se^2
(\C^k) \oplus (\C^k \otimes \C^{m-k})^* \oplus \Se^2 (\C^{m-k}).$
Hence, by Tables IIa and IIb we have $k=m-k=1$ which implies $\dim
\mf l=2< \dim \Se^2 (\C^2).$ Summing up we have the following
minimal subalgebra: $ \un 1$ acting on  $\Sp (1) / \U (1)$ and
$\su m$ acting on  $\Sp (m) / \U (m).$
%
%
\subsection{The case $K= \SO (p) \otimes \Sp (q),\ pq=m,\ p\geq 3,\ q \geq 1$}
The dimension of a Borel subgroup of $K^{\C}$ is equal or less
than $\frac{p^2}{4} + q^2 + q,$ while $\dim M =\frac{1}{2}
pq(pq+1),$ since $m=pq.$ Now, let $f(x)=x^2(2q^2 -1) +2xq
-4q^2-4q.$ Then $f'(x)>0$ for $x \geq 0$ and $f(3)>0$ since $q\geq
1.$ Then the $K-$action cannot be coisotropic.
%
%
\subsection{The case $K=\Sp (k) \times \Sp (m-k)$} \label{cc}
Since $K$ is a symmetric subgroup of $\Sp(m),$ the $K-$action is
hyperpolar. We shall analyze the subgroups of $K.$ The manifold
$M$ parametrizes the space of Lagrangian subspaces of $\C^{2m}$
respect to a symplectic form. We consider $ \omega (X,Y)= X^t J Y
$ where
$$
J= \left(
\begin{array}{c|c|c|c}
0   & -\I_k  & 0   & 0  \\\hline \I_k & 0     & 0   & 0
 \\\hline
0    & 0     & 0   & -\I_{m-k} \\\hline
0    & 0     & \I_{m-k} & 0   \\
\end{array}
\right) = \left(
\begin{array}{c|c}
J_k   & 0   \\\hline
0 & J_{m-k}  \\
\end{array}
\right)
$$
Let $W_o =< e_1, \ldots , e_k > \oplus <e_{m+k+1}, \ldots,
e_{2m}>.$ Notice that $W_o$ is a Lagrangian subspace of $\C^{2m},$
$< e_1, \ldots , e_k >$ ($<e_{n+k+1}, \ldots, e_{2n}>$) is a
Lagrangian subspace of $\C^{2k}$ ($\C^{2(m-k)}$ ) respect to the
symplectic form $\omega_k = X^t J_k  Y$ ($\omega_{m-k} (X,Y)= X^t
J_{m-k} Y$). Hence the orbit of $K$ through $W_o$ is $\Sp(k) / \U
(k) \times \Sp (m-k) / \U(m-k),$ and the tangent space at $[\U
(m)]$ splits
$$
\Se^2 (\C^m)= \Se^2 (\C^{k} ) \oplus \Se^2 (\C^{m-k}) \oplus (
\C^{k} \otimes \C^{m-k})^*,
$$
as $\U(k) \times \U(m-k)-$modules, proving that the slice
representation is given by $\C^{k} \otimes \C^{m-k}$ on which
$\U(k) \otimes \U(m-k)$ acts. Note that the slice appears in Table
Ia: this is another way to prove that the $K-$action is
multiplicity free. Now let $L\subseteq K= \Sp (k) \times \Sp(m-k)$
and let $\mf l$  be the Lie algebra of $L.$ Suppose $\mathfrak{l}$
acts coisotropically. We consider the projections $\sigma_1 :\mf l
\longrightarrow \mf{sp} (k),$ $\sigma_2 :\mf l \longrightarrow
\mf{sp} (m-k)$ and we put $\mf{l}_i= \sigma_i (\mf l ).$ This
means that $\mf l \subset \li1 + \li2,$ $\li 1 + \li 2$ acts
coisotropically on $\Sp (m) / \U (m),$ so $\li 1,$ respectively
$\li 2,$ acts coisotropically on $\Sp (k) / \U (k),$ respectively
$\Sp (m-k) / \U (m-k).$ Then we have the following
possibility \\
$\ $ \n
{\bf  \S 1 $\li 1$ and $\li 2$ act both transitively} \\
Hence $\mf l = \asp(k) + \asp (m-k)$ or $\mf l = \asp (k) + \theta
(\asp (k)),$ where $\theta$ is an automorphism of $\asp (k).$ The
first case corresponds to $\Sp (k) \times \Sp(m-k)$ that we have
just considerated. The second case must be excluded by dimensional
condition. Indeed, the dimension of a Borel subgroup of $\mf
l^{\C}$ is $k^2+k$ while $ \dim \Sp(2k)/ \U(2k)=2k^2 + k $
$\ $ \\
{\bf  \S 2 $\li 1$ acts transitively and $\li 2$ acts coisotropically} \\
We must consider the following cases
\begin{enumerate}
\item $\li 1= \asp (k)$ and  $\li 2$ has a fixed point. Hence $\mf
l = \li 1 + \li 2.$ The orbit through $W_o$ is a complex orbit and
the slice is given by $\Se^2 (\C^{m-k}) \oplus (\C^{m-k} \oplus
\C^k)^*$ on which $\un k$ acts on $\C^k$ and $\li 2$ acts on
$\C^{m-k}.$ By Tables IIa and IIb, this representations fails to
be multiplicity free when $m-k \geq 2$. If $m-k=1,$ note that $\li
2$ must be $\un 1,$ then the action is multiplicity free but the
scalar cannot be removed. Summing up, we have the following
multiplicity free action: $\mf l=\asp (m-1) + \un 1$
\item $\li 2 \subseteq \asp (m_1) + \asp(m_2),$ where $m_1 +
m_2=m-k.$ We may suppose, up to conjugation in $\asp(m)$, $ k \geq
m_1 \geq m_2.$ Let $\li 2= \asp (m_1) + \asp(m_2).$ Then $\mf l=
\li 1 + \li 2,$ which corresponds to $L=\Sp(k) \times \Sp (m_1)
\times \Sp(m_2) \subseteq K= \Sp (k) \times \Sp(m-k).$ We have
proved that there exists $W \in \Sp (m-k) / \U (m-k)$ such that
$\Sp (m_1) \times \Sp(m_2) W$ is a complex orbit. Since $\Sp(k)
\times \Sp(m-k) W_o= \Sp(k) / \U(k) \times \Sp(m-k) / \U(m-k),$
the orbit $\Sp(k) \times \Sp (m_1) \times \Sp(m_2)W$ is a complex
orbit and the slice is given by
$$
(\C^k \otimes \C^{m_1})^* \oplus (\C^k \otimes \C^{m_2})^* \oplus
(\C^{m_1} \otimes \C^{m_2})^*
$$
on which $\U(k)$ acts on $\C^k,$ $\U(m_1)$ acts on $\C^{m_1}$ and
$\U (m_2)$ acts on $\C^{m_2}.$ By Tables IIa and IIb we must
assume $m_1=m_2=1,$ so the slice becomes $(\C^k \oplus \C^k \oplus
\C)^*$ and the two copies of $\U(1)$ act as $(e^{-i \psi}, 1 ,
e^{-i \psi})$ and $(1,e^{-i \phi}, e^{-i \phi})$  respectively.
Since a representation $(\rho,V)$ is multiplicity free if and only
if the dual representation ($\rho^*, V^*)$ is, we may assume that
$S=\C^k \oplus \C^k \oplus \C.$ To solve this case we apply $(ii)$
of Theorem \ref{class}. Note also by the Theorem 1.1 page $7$ in
\cite{Kol} we may analyze the slice representation.  Firstly, let
$ 1 \in \C.$ The orbit is ${\rm S}^1$ and the slice is given by
$\R \oplus \C^k \oplus \C^k$ on which $\U(1) \times \U(k)$ acts as
follows: $(e^{i \phi}, A) (\alpha,v,w)=(\alpha, e^{i \phi}Av,e^{-i
\phi}Aw).$ Now, we consider $(0,0,(1, \ldots,0));$ the orbit is
the unit sphere and the slice becomes $\R \oplus \R \oplus \C
\oplus \C^{k-1}$ on which $\T^1 \times \U(k-1)$ acts as follow:
$(e^{i \phi}, A) (\alpha, \beta, z, v)=(\alpha, \beta, e^{i
\phi}z, Av).$ Now it is easy to see that $H_{{\rm princ}}
=\U(k-2)$ and the cohomogeneity is $4,$ thus proving
$$
4=\mathrm{ch}(H, S)= \mathrm{rank}(H) - \mathrm{rank}(H_{\rm
princ} ) =2+k-(k-2).
$$
We must analyze the behaviour of the subgroup of $H.$ However, by
the Restriction lemma \cite{Hw}, if one takes
 $L \subset \Sp(1) \times \Sp(1)$ such that
$\Sp(m-2) \times L$ acts coisotropically on $\Sp(m)/ \U(m)$ then
$L$ acts coisotropically on $\Sp(2) / \U(2).$ Hence, for
dimensional reasons, $L$ must be $\U(1) \times \Sp(1)$. However,
the orbit through $W$ is a complex orbit and the slice becomes $
(\C^k \oplus  \C^k \oplus \C^k \oplus \C^k)^* \oplus (\C)^* \oplus
(\C)^* \oplus \Se^2 (\C) $ on which $\U(k)$ acts on $\C^k,$ so by
Tables IIa and IIb the action fails to be multiplicity free.
\end{enumerate}
{\bf  \S 3 $\li 1$ and $\li 2$ act both coisotropically } \\
Since if both $\li 1$ and $\li 2$ have a fixed point, then $\mf l
\subseteq \li 1 + \li 2$ has a fixed point, we shall analyze the
following cases: $\li 1=\un k,$ $\li 2=\asp(m_1) + \asp(m_2)$ and
$\li 1= \asp(k_1)+ \asp(k_2),$ $\li 2= \asp(m_1) + \asp(m_2).$
Since $\li 1 + \li 2 \subseteq \asp(k) + \asp(m_1) + \asp(m_2),$
we have $m_1=m_2=1.$ In particular, the first case must be
excluded for dimensional reason. In the second case $\mf l= \li 1
+ \li 2,$ which corresponds to $L=\Sp(k_1) \times \Sp(k_2) \times
\Sp(1)\times \Sp(1)$ and one may prove that $L$ has a complex
orbit given by $\Sp(k_1) / \U(k_1) \times \Sp(k_2)/ \U(k_2) \times
\Sp (1) / \U(1) \times \Sp(1)/ \U(1)$ whose slice representation
fails to be multiplicity free.
%
%
%
%
%
%
\section{$M=\SO (2m)/ \U (m)$}
In the following subsections we will go through all maximal
subgroups $K$ of $\SO (2m)$ according to Table {\bf IV} in the
Appendix.
%
%
\subsection{The case $K =\rho (H),$ $H$ simple such that $\rho \in Irr_{\R}(H),$
$d(\rho)=2m$}

By lemma \ref{zio} we shall analyze the cases $(\mf{so} (6),
\Lambda_3)$ and $(\mf{so} (6), \Lambda_2),$ which correspond to
the transitive  action of $\SO(6)$ on $\SO (6) / \U(4).$
%
%
%
\subsection{The fixed point case $K= \U(m)$}
We use the same notation and the same strategy as in section
\ref{fisso} By Table Ia $\U (m)$ acts coisotropically on
$\Lambda^2 (\C^m )$ and the scalar can be reduced. Throughout this
section we denote by $\mf z$ the center of Lie($\U(m)$)$=\un m$
and by $\mf t_m$ the Lie algebra of a maximal torus of $\su m
\subseteq \un m.$ Let $L$ be a compact subgroup of $\U(m)$ such
that $\mf l= \mf z + \li 1,$ where $\li 1$ is a maximal Lie
algebra of $\su m$ (see Table {\bf V} in the Appendix). By lemma
\ref{zio} the case $\li 1= \so m$ can be excluded, while the case
$\li 1 = \asp (n), 2n=m,$ appears when $n=2$ and the slice becomes
$\C \oplus \C^5$ on which $\Sp(2) / \Z_2= \SO (5)$ acts on $\C^5.$
Then $\mf l=\mf z + \asp (2)$ acts coisotropically and the scalar
cannot be removed. Notice that, since the slice of the orbit
through $1 \in \C$ is $\R \oplus \C^5$ on which $\SO(5)$ acts on
$\C^5,$ one may prove, see also \cite{He}, the slice fails to be
polar. This case is maximal, since for every  $\mf h \subseteq \mf
\asp (2)$ we have $\mf z + \mf h$ does not satisfy the dimensional
condition.

If $\li 1= \R + \su k + \su{m-k}$ then the slice becomes $
\Lambda^2 (\C^m)= \Lambda^2 (\C^{m-k}) \oplus (\C^k \oplus
\C^{m-k} )^* \oplus \Lambda^2 (\C^{m-k} ), $ on which $\su k ,$
respectively $\su{m-k},$ acts on $\C^k,$ respectively $\C^{m-k}.$
Hence by Tables Ia, Ib and Tables IIa, IIb, we have $k=1$ and the
slice becomes $\Lambda^2 (\C^{m-1}) \oplus (\C \otimes \C^{m-1}
)^* .$ The scalars, $\mf z$ and $\R=\mf a,$ the centralizer of
$\su{m-1}$ in $\su m \subseteq \un m,$ act as follows: let $(\psi,
\theta) \in \mf a \times \mf z$, then $(\psi, \theta) (v,w)= (e^{
2 i( \theta  - \frac{1}{m-1} \psi)} v, e^{ -i (2 \theta +
\frac{m-2}{m-1} \psi) } w).$ Hence, the action is multiplicity
free and we shall show how many centers we need. Firstly, we
assume $m \geq 5.$ By Table IIa the scalars can be reduced in the
following cases: when $m-1$ is even, we need only a one
dimensional center acting on the first submodule, that is
satisfied with the line $\R( \frac{1}{m-1} ),$ where $\R(\alpha)$
means every line in the plane $(x,y) \in \mf a \times \mf z$
different from $y=\alpha x,$ while when $m-1=2s+1$ one may prove
that we can reduce the scalars, but the scalars cannot be removed.
When $m=4,$ the slice becomes $( \C^3 \oplus \C^3 )^*,$ so by
Table IIa, the scalars cannot be removed, but can be reduced if
the center acts as $(z^a, z^b)$ with $a \neq b.$ This corresponds
to $\R(0) + \su 3 .$ Finally, when $m=3,$ the slice becomes $\C
\oplus \C^2$ and it is easy to see that the minimal subalgebra is
$\mf z + \su 2.$ Notice that for $m \geq 4$ these actions are
maximal by Tables IIa and IIb. If $m=3,$ then also $\mf z + \mf
t_3$ acts coisotropically on $\SO(6) / \U(3)$ and when $m=2$ we
have also $\R(0),$ line in $\mf a \times \mf z,$ acting on $\SO(4)
/ \U(2).$

The case (iv) can be excluded by dimensional condition as in
section \ref{fisso}. Indeed,
$\dim_{\C}\SO(2m)/\U(m)=\frac{1}{2}pq(pq-1),$ since $m=pq$ and the
dimension of a Borel subgroup of $(\SU(p) \otimes \SU(q))^{\C}$ is
$\frac{1}{2}(p^2 +q^2 +p +q-4)$. We shall prove that $pq(pq-1)
> p^2+q^2+q+p-2$ which implies that the dimensional condition is not
satisfied for a Lie group with Lie algebra $\mf z + \su p + \su
q$.

Let $f(x)=x^2(q^2-1)-x(q+1) -q^2 -q+2.$ Then $f'(x)=2x(q^2-1)
-q+4>0,$ for $x \geq 3$ and $f(3)=9(q^2-1)-3(q+1)-q^2-q+4>0,$ when
$q \geq 2.$

Finally, we consider the case (v). By lemma \ref{zio} we have only
the case $\su m$ which has just been analyzed. Summing up, if $L
\subset \U (m)$ acts coisotropically on $M$ then, up to
conjugation in $\mf{o} (2m),$ the minimal algebra are in the
following table
\begin{small}
\begin{table}[ht]
\renewcommand\arraystretch{1.0}
\noindent\[
\begin{array}{|l|c|l|}
\hline
 \mf l      &  M   & {\rm conditions} \\ \hline\hline
\R (0)      & \SO (4) / \U(2) & \R(0) \ \mathrm{line\ in} \ \mf
t_2 \times z
                                     \\ \hline
\mf z + \mf t_3 & \SO (6) / \U(3) &  \\ \hline
\R (\frac{1}{2k} ) + \su{2k}   & \SO (4k+2) / \U (2k+1)  &  k\geq
2,\ \R(\frac{1}{2k}) \ \mathrm{line \ in} \ \mf a \times \mf z \\
\hline \R  + \su{2k+1}   & \SO (4k+4) / \U (2k+2)  & k\geq 2,\ \R
\ \mathrm{means \ any \ line \ in} \ \mf a \times \mf z \\ \hline
\R (0) + \su 3 & \SO (8) / \U (4) & \R(0) \ \mathrm{line \ in} \
\mf a \times \mf z
\\ \hline
\mf z + \su 2   & \SO (6) / \U (3) & \\ \hline
\su m       & \SO(m) / \U(m) & m \geq 2 \\ \hline
\mf z + \asp (2)  & \SO (8) / \U (4)  &  \\ \hline
\end{array}
\]
\end{table}
\end{small}
%
%
%
\subsection{The case $K=\SO (p ) \otimes \SO (q),\ 3 \leq p \leq q$}
By a straitforward calculation one may prove that $\SO (p) \otimes
\SO (q),\ 3 \leq p \leq q $ does not satisfy  the  dimensional
condition.
%
%
%
\subsection{The case $K= \Sp(p) \otimes \Sp(q),\ 4pq=2m$ }
One may prove that $K$ does not satisfy the dimensional condition
unless $p=q=1$ and $p=1$ and $q=2.$ Now, the case $\Sp(1) \otimes
\Sp(1)$ corresponds to the transitive action of $\SO(4)$ on
$\SO(4)/ \U(2),$ while $\Sp(1) \otimes \Sp(2)$ acts on $\SO(8) /
\U(4).$ Since $\Sp(1) \otimes \Sp(2) \cap \U (4)= \T^1 \cdot \Sp
(2)$ the $\Sp(1) \otimes \Sp (2)-$orbit through $[\U (4)]$ is a
complex orbit and the slice is given by $ \C^5, $ on which
$\Sp(2)$ acts on $\C^5$  as $\Spin (5) / \mathbb{Z}_2= \SO (5).$
By Table Ia the action is multiplicity free and the scalar cannot
be removed. Thanks to dimensional condition we must analyze only
the following subgroups of $\Sp(1) \otimes \Sp(2):$ $H=\T^1 \times
\Sp (2),$ which has been considerated in the fixed point case, and
$H=\Sp(1) \otimes (\Sp(1) \times \Sp(1)).$ However, $H \cap \U(4)=
\T^1 \cdot (\Sp (1) \times \Sp(1))$ and the slice becomes
$\Lambda^2 (\C^2) \oplus (\C^2 \otimes \C^2)^* \oplus \Lambda^2
(\C^2)$ on which $\Sp(1) \otimes \Sp(1)$ acts on $\C^2 \otimes
\C^2.$ By Table IIb, we need two dimensional scalars acting on
$\C^2 \otimes \C^2,$ hence the action fails to be multiplicity
free.
%
%
\subsection{The case $K= \SO (k) \times \SO(2m-k)$}
Since $K$ is a symmetric group of $\SO (2m),$ the $K$-action is
hyperpolar. We shall analyze the behaviour of the closed subgroups
of $K=\SO (k) \times \SO(2m-k),$ so it is very useful to get a
complex orbit of $K.$ Notice that we may assume $k\leq m.$
Firstly, we suppose $k=2s.$ The homogeneous space $M=\SO(2m) /
\U(m)$ parametrizes the almost complex structure $\R^{2m}$ that
are orthogonal  and compatible with a fixed orientation. Let
$J_1,$ respectively $J_2,$ be almost complex structure of
$\R^{2s},$ respectively $\R^{2(m-s)},$ as above and let $\J_o =J_1
\oplus J_2.$ Clearly, $\J_o$ is an orthogonal almost complex
structure of $\R^{2m},$ the orbit $K \J_o$ is $ \SO(2s) / \U (s)
\times \SO (2m-2s) / \U (m-s) $ and the slice is given by $ \C^s
\otimes \C^{m-s} $ on which $\U(s)$ acts on $\C^s$ and $\U(m-s)$
acts on $\C^{m-s},$ i.e. $K \J_o$ is a complex orbit. If $k=2s+1$
we split $\R^{2n}= \R^{2s} \oplus \R^2 \oplus \R^{2(m-s-1)}$ and
we consider $\J_e=J_1 \oplus J_2 \oplus J_3,$ where $J_1,$ $J_2$
and $J_3$ are orthogonal almost complex structures of $\R^{2s},$
$\R^2$ and $\R^{2(m-s-1)}$ respectively. One may prove that the
orbit through $\J_e$ is $\SO(2s+1) / \U (s) \times \SO (2(m-s-1)
+1) / \U (m-s-1) ,$ and the slice is given by $(\C^s \otimes
\C^{m-s-1})^*.$

Now let $L\subseteq K= \SO (2s) \times \SO (2m-2s)$ and let $\mf
l$  be the Lie algebra of $L.$ Suppose $\mathfrak{l}$ acts
coisotropically. We consider the projections $\sigma_1 :\mf l
\longrightarrow \mf{so} (k),$ $\sigma_2 :\mf l \longrightarrow
\mf{so} (2m-k)$ and we put $\mf{l}_i= \sigma_i (\mf l ).$ This
means that $\mf l \subset \li1 + \li2,$ $\li 1 + \li 2,$  acts
coisotropically on $\SO (m) / \U (m),$ so $\li 1,$ respectively
$\li 2,$ acts coisotropically on $\SO (2s) / \U (s),$ respectively
on $\SO (2m-2s) / \U (m-s).$ In the sequel we refer to Tables Ia,
Ib and Tables IIa, IIb in the Appendix, for all the conditions
under which one can remove or reduce the scalar preserving the
multiplicity free action. Then we have the following possibility
\\
%
%
\noindent
{\S \bf 1 $\li 1$ and $\li 2$ act both transitively} \\
For dimensional reasons $\mf l = \mf{so}(2s) + \mf{so}(2(m-s))$
which has just
been considerated.\\
%
%
\noindent
{\S \bf 2 $\li 1$ acts transitively and $\li 2$ acts coisotropically} \\
We must analyze the following cases
\begin{enumerate}
\item $\li 1 = \mf so (2m-2)$ and $\li 2=0 \subseteq \so 2.$ The
orbit through $\J_o$ is complex and the slice becomes $ (\C
\otimes \C^{m-1} )^* $ where $\un{m-1}$ acts on $\C^{m-1}.$ Hence,
by Table Ia, the action is multiplicity free. Since the
cohomogeneity is $1$ this action is hyperpolar. \item $\li 2$ has
a fixed point. The orbit through $\J_o$ is a complex orbit $\SO
(2s) / \U (s),$ so we are going to analyze the slice
representation according the table appears in section 4.2.
\begin{itemize}
\item $\li 1= \R (0) \subseteq \un 2 \subseteq \mf{so}(4).$ The
slice becomes
$$
(\C^s \otimes \C)^* \oplus (\C^s \otimes \C)^* \oplus (\C \otimes
\C)^* .
$$
on which $\un s$ acts on $\C^s$ and $\R (0)$ acts on $\C .$ Hence
the action fails to be  multiplicity free since the scalars act on
$(\C^s \otimes \C)^* \oplus (\C^s \otimes \C)^*$ as a one
dimensional scalar; \item the cases $\li 2 = \mf z +\mf t_3,$ $\li
2 = \R(\frac{1}{2k})  +\su{2k},$ $\li 2 = \R + \su{2k+1},\ k \geq
2.$ $\li 2 = \R (0) + \su 3$ and $\li 2 = \mf z +\asp (2)$ can be
excluded since two many terms appear in the slice. Indeed, for
example, let  $\li 2= \R(\frac{1}{2k})  +\su{2k}.$ Then $\mf l=
\li 1 + \li 2,$ and the slice becomes $(\C^s \oplus \C^{2k})^*
\oplus (\C^s \oplus \C)^* \oplus \Lambda^2 (\C^{2k}) \oplus (\C
\oplus \C^{2k})^* .$ By Tables IIa and IIb this action is not
multiplicity free. \item $\li 2 \subset \mf z + \su{m-s}.$ The
slice becomes $(\C^{s} \otimes \C^{m-s})^* \oplus \Lambda^2
(\C^{m-s})$ on which $\un s$ acts on $\C^s$ and $\li 2$ acts on
$\C^{m-s}.$ If $m-s\geq 4$ then the action fails to be
multiplicity free while if $m-s=3$ or  $m-s=2$ then the action is
multiplicity free with the scalar $\mf z.$ Summing up we have the
following subalgebra: $\so{2m-6} + \un 3,\ m \geq 5,$ and
$\so{2m-4} + \un 2, \ m \geq 4$ acting on  $\SO(2m) / \U(m).$
\end{itemize}
%
%
%
\item $\li 2= \asp(1) + \asp(2).$ Then $\mf l = \li 1 + \li 2$ and
a complex orbit is given by $ \SO(2(m-4) / \U(m)) \times \C.$
However, one may prove that the slice fails to be multiplicity
free; \item $\li 2 \subseteq \so{m_1} + \so{m_2}.$ We may assume,
up to conjugation, that $2s \geq m_1 \geq m_2.$ Let $\li 2=
\so{m_1} + \so{m_2}.$ Then $\mf l= \li 1 + \li 2$ which
corresponds to $\SO(2s) \times \SO(m_1) \times \SO(m_2 ).$ Assume
both $m_1$ and $m_2$ are even. We know that there exists $\J_o$
such that $\SO (m_1) \times \SO(m_2) \J_o$ is a complex orbit in
$\SO(2m-2s) / \U(m-s).$ Hence $\SO (2s) \times \SO(m_1) \times
\SO(m_2) \J_o$ is a complex orbit and the slice is given by
$$
( \C^s \otimes \C^{\frac{m_2 -1}{2}} )^* \oplus (\C^s \otimes
\C^{\frac{m_1-1}{2}})^* \oplus ( \C^{\frac{m_1-1}{2}} \otimes
\C^{\frac{m_2-1}{2}})^*.
$$
Since $s\geq 2,$ by Tables IIa and IIb we get $m_1=m_2=2$ and the
slice becomes
$$
(\C^s \otimes \C)^* \oplus (\C^s \otimes \C)^* \oplus (\C \otimes
\C)^*
$$
on which $\U (s)$ acts on $\C^s.$ The center of $\U (s)$ acts as
as $(e^{-i \theta}, e^{-i \theta}, 1),$ while $\SO(2) \times
\SO(2)$ acts as $(e^{-i \phi}, e^{-i \psi}, e^{-i (\phi+ \psi)}).$
Hence, we get the following minimal subalgebra: $\so 4 + \R + \R$
acting on $\SO(8) / \U(4)$ and $\so{2s} + \R(1,-1),$ where
$\R(1,-1)$ is a line different form $y=x, y=-x,$ acting on
$\SO(2(s+2))/ \U(s+2),$ for $s \geq 3.$

Finally, assume that $m_1$ and $m_2$ are odd. Notice that the case
$m_1=m_2=1$ has been considerated. Hence the slice of the complex
orbit $\SO(2s) \times \SO(m_1) \times \SO(m_2) \J_e$ is given by $
(\C^s \otimes \C^{ \frac{m_1-1}{2} } )^* \oplus (\C^s \otimes \C^{
\frac{m_2-1}{2} } )^* \oplus (\C^s \oplus \C )^* \oplus ( \C^{
\frac{m_1-1}{2} } \otimes  \C^{ \frac{m_2-1}{2} } )^* $ so this
action is not multiplicity free.
\end{enumerate}
%
%
%
\noindent
{\S \bf 3 $\li 1$ and $\li 2$ act both coisotropically} \\
As in section \ref{cc} we may prove that $\mf l$ does not act
coisotropically. For example, let $\li 1 = \un l $ and let $\li 2
= \mf{so}(p) + \mf{so}(q),$ where $p,q$ are even. Then $\mf l =
\li 1 + \li 2$ and the orbit through through $\J_o' \oplus \J_o$
is a complex orbit whose slice is given by $\Lambda^2 (\C^l)
\oplus (\C^{\frac{p}{2}}  \otimes \C^{ \frac{q}{2} } )^* \oplus
(\C^{l} \otimes \C^{\frac{p}{2}}  )^* \oplus,
 (\C^{l} \otimes \C^{\frac{q}{2}}  )^*,$
on which $\un l$ acts on $\C^l,$ and $\un{\frac{p}{2}},$
respectively $\un{\frac{q}{2}},$ acts on $\C^{\frac{p}{2}},$
respectively $\C^{\frac{q}{2}}.$ Hence, this action fails to be
multiplicity free.

Now we are going to analyze the behaviour of the subgroup of
$\SO(k) \times \SO (2m-k)$ when $k$ is odd. The maximal subgroup
$L$ of $\SO(k) \times \SO (2m-k)$ are: $H \times \SO(2m-k),$ where
$H$ is a maximal subgroup of $\SO(k),$ $\SO(k) \times H$ where $H$
is a maximal subgroup of $\SO(2m-k)$ and when $k=2m-k,$ $\SO(k)
\times A(\SO(k))$ where $A$ is an automorphism of $\SO(k).$
However, the last case can be excluded for dimensional reasons.

Since $k$ is even we have the following cases: $H= \SO(p) \otimes
\SO(q),$ $pq=k,$ $3 \leq p \leq q$ and $H =\sigma (L),$ $L$ simple
such that $\sigma \in Irr_{\R}(L).$ The first case may excluded by
dimensional reasons. Indeed, if $H \times \SO(2m-k)$ acts
coisotropically on $M= \SO(2m) / \U(m)$ then, by Restriction
lemma, see \cite{Hw}, $H \times \SO(2m-k)$ acts coisotropically on
the complex orbit of $\SO(k) \times \SO(2m-k),$ that is $\SO(2s+1)
/ \U(s) \times \SO(2(m-s-1)+1) / \U(m-s-1),$ since $k=2s+1.$ In
particular $H$ acts coisotropically on $\SO(2s+1) / \U(s).$
However the dimension of a Borel subgroup of $H^{\C}$ is lesser
than $\frac{p^2 + q^2}{4}$ while $\dim \SO(2s+1) /
\U(s)=\frac{p^2q^2-1}{8},$ since $s=\frac{pq-1}{2}.$ The
inequality $ 2(p^2 +q^2) < p^2q^2 -1 $ means that the dimensional
condition does not satisfy.

Let $f(x)=x^2(p^2-2) -2p^2 -1.$ Then $f'(x) >0$ if $x>0$ and
$f(3)=p^2-19>0.$ Hence the action fails to be multiplicity free.

Now, we shall prove that if $H =\sigma (L),$ $L$ simple such that
$\sigma \in Irr_{\R}(L)$ then $H=G_2 \subseteq \SO(7).$ As before,
if $H \times \SO(2m-k)$ acts coisotropically then the dimension of
a Borel subalgebra of $\mf h$ must satisfy the following
inequality
\begin{equation}\label{eccolo}
\dim \mf b \geq \frac{d^2-1}{8}
\end{equation}
We may analyze any simple Lie algebra as in lemma \ref{zio}.
Notice that $d(\sigma)$ must be odd. This is a straitforward
calculation and easy to check. We demonstrate our method analyzing
the cases $\mf h= \su m$ and $\mf g_2.$

If $\mf h = \su m,$ then $\dim \mf b = \frac{1}{2}(m-1)(m+2).$ The
case $m=2$ give rise a real representation $2 \Lambda_1$ which
corresponds to the transitive action of $\SO(3).$
Now, assume $m\geq 3.$ It is well known that if $\sigma=
\sum_{i=1}^{m-1} a_i \Lambda_i$ is a contragradient representation
then $a_i=a_{m-i},$ and one may prove that $d(\sigma) \geq
d(\Lambda_1 + \Lambda_{m-1}).$ Since $d( \Lambda_1 +
\Lambda_{m-1})=m^2 -1 \geq \frac{5}{2}m,$ (\ref{eccolo}) does not
hold for any real representation.
%
Assume $\mf h=\mf g_2.$ Since the dimension of a Borel subalgebra
is $8$ hence $(1)$ becomes $63 \geq d^2 (\sigma)$ that is verified
only for $\Lambda_1$ which corresponds to ${\rm G}_2 \subseteq
\SO(7)$ acting on $M=\SO(8)/ \U(4).$ Since ${\rm G}_2 \cap \U
(4)=\SU (3),$ the orbit through $[\U (4)],$ ${\rm G}_2 / \SU (3)
\cong S^6,$ is totally real. Indeed, let $\phi: \SO(8) / \U(4)
\lra \mf{g}_2^*$ be the moment map. Then ${\rm G}_2
\phi([\U(4)])={\rm G}_2 / P$ is a flag manifold, and $\SU (3)
\subseteq P.$ However $\SU (3)$ is a maximal subgroup of ${\rm
G}_2$ so $P=G_2$ and $\phi([\U (4)])=0.$ Now, it is easy to check
that ${\rm G}_2 [\U (4)]$ is totally real. Moreover, since $2
\dim_{\R} {\rm G}_2 / \SU(3)= \dim_{\R} \SO(8)/ \U(4),$ the slice
representation can be deduced immediately from the isotropic
representation of $\SU(3)$ on ${\rm G}_2 / \SU(3),$ showing that
the cohomogenity of the $G_2 -$action is $1,$ which implies ${\rm
G}_2$ acts hyperpolarly on $\SO(8)/ \U(4).$

Now shall investigate ${\rm G}_2 \times \SO(2s+1),$ for every $s
\geq 1,$ acting on $\SO(2(s+4)) / \U(s+4).$ The isotropy group of
${\rm G}_2 \times \SO(2s+1) \J_o,$ is $\SU (3) \times \U(s)$ and
the slice, from real point of view, is given by $ \C^3 \oplus
(\C^3 \otimes \C^s) $ on which $\SU(3)$ acts on $\C^3$ and $\U(s)$
acts on $\C^s.$ We shall prove that (ii) of Theorem \ref{class} is
not satisfied. By the slice theorem, see \cite{Kol}, it is enough
to study the slice representation.

The case $s=1$ is a straitforward calculation and for dimensional
reasons we shall assume $s \geq 3.$ Let $v\in \C^3$ and let $w \in
\C^s$ be two unit vectors. One can prove that the isotropy group
of $v + v \otimes w$ is $\SU(2) \times \U(s-1)$ which acts on the
slice $\C^2 \oplus \C^2 \otimes \C^{s-1}.$ If we iterate this
procedure two times then we get that the regular isotropy is
$\U(s-3)$ and the cohomogeneity is $7.$ However $7 \neq {\rm
rank}( {\rm G}_2 \times \SO(2s+1) )- {\rm rank}(\U(s-3))=5$.
%
%

Finally, we shall analyze ${\rm G}_2 \times {\rm G}_2,$ acting on
$\SO(14)/ \U(7).$ However, for dimensional reasons, the action
fails to be multiplicity free.
%
%
%
%
%
\section{$M=\Ea_7 / \T^1 \cdot \Ea_6$}
In this section we analyze the behaviour of the subgroup of
$\Ea_7.$ By dimensional condition, a subgroup $K\subseteq \Ea_7$
which acts coisotropically on $M$ must satisfy $\dim K \geq 47.$
The maximal subgroups of $\Ea_7$ which satisfy the above
inequality (see \cite{Kol} page 41) are the following
$$
\begin{array}{|l|c|c|c|} \hline
{\rm maximal \ rank} & \T^1 \cdot \Ea_6
                     & \SU(2) \cdot \Spin (12)
                     & \SU(8) / \Z_2  \\ \hline
{\rm no \ maximal \ rank} & \SU(2) \cdot \Fa_4  & & \\ \hline
\end{array}
$$
We are going to analyze these cases separately.
%
%
\subsection{The fixed point case $K=\T^1 \cdot \Ea_6$}
The subgroup $K$ acts coisotropically, since it has a fixed point
and the slice representation, which is given by $(\C^{27},
\Lambda_1),$ appears in Table Ia. Note also that the scalar cannot
be removed. The unique maximal subgroup $H$ of $\T^1 \cdot \Ea_6 $
which satisfies $\dim H \geq 47 $ is $\T^1 \cdot \Fa_4 .$ However
this actions fails to be multiplicity free. Indeed, the slice
representation is given by $\C^{26} \oplus \C,$ (see \cite{Ad}
lemma 14.4 page 95) so by Table Ia this actions fails to be
multiplicity free.
%
%
\subsection{The case $K=\SU(2) \cdot \Fa_4$}
By Table 25 in \cite{di2} page 204, one sees, after conjugation,
$\Fa_4$ is contained in $\Ea_6.$ Hence the connected component of
$K \cap \T^1 \cdot \Ea_6 $ is $\Fa_4$ or $\T^1 \cdot \Fa_4$, since
$K$ is a maximal subgroup. However $\C^{27}=\C^{26} \oplus \C$ as
$\Fa_4 -$submodules (see Lemma 14.4 page 95 \cite{Ad}). Hence
 $K \cap \T^1 \cdot \Ea_6 = \T^1 \cdot \Fa_4$ so the orbit through
$[{\rm T}^1 \cdot \Ea_6]$ is a complex orbit which slice
representation fails to be multiplicity free.
%
%
\subsection{The case $K=\SU(2) \cdot \Spin (12) $}
$K$ is a symmetric group of $\Ea_7$ hence the action is hyperpolar
on $M.$ Now, since any automorphism of $\Ea_7$ is an inner
automorphism then for any $\sigma, \tau \in \Aut (\Ea_7)$ there
exists an element $g \in \Ea_7$ such that $\sigma$ and $Ad(g^{-1})
\circ \tau \circ Ag(g)$ commute. Hence we may assume that $K \cap
\T^1 \cdot \Ea_6$ is a symmetric subgroup of $K$ and $\T^1 \cdot
\Ea_6.$ Since the symmetric subgroup of $\Ea_6$ are the following
$$
\begin{array}{|c|c|c|c|} \hline
 \T^1 \cdot \Spin (10) &
 \T^1 \cdot \SU (6)    &
 \Fa_4                 &
 \Sp (4) / \Z_2        \\ \hline
\end{array}
$$
then $K \cap \T^1 \cdot \Ea_6 = \T^1 \cdot ( \T^1 \cdot
\Spin(10)),$ where the first $\T^1$ lies in $\SU(2),$ but it is
different from the centralizer of $\Ea_6$ in  $\Ea_7,$ while the
second is the centralizer of $\Spin(10)$ in $\Spin(12).$ The slice
representation is given by $\C^{16}$ on which $\T^1 \cdot \T^1
\cdot \Spin (10)$ acts. Hence $K$ acts coisotropically on $M$.

Now we analyze the behaviour of the subgroup of $K.$

Let $L=\T^1 \cdot \Spin (12),$ where $\T^1\subseteq \SU(2).$ Then
$\T^1 \cdot \Spin (12) \cap \T^1 \cdot \Ea_6 =\T^1 \cdot (\T^1
\cdot \Spin (10))$ and the slice becomes $\C^{16} \oplus \C,$ on
which $\T^1 \cdot (\T^1 \cdot \Spin(10))$ acts. Note that the
first scalar acts on $\C$ while the centralizer of $\Spin(10)$ in
$\Spin(12)$ does not. Hence, the action is multiplicity free,
since the $\Spin(10)-$action on $\C^{16}$ is multiplicity free.

The case $L=\Spin(12)$ must be excluded, since $L \cap \T^1 \cdot
\Ea_6= \T^1 \cdot \C^{16},$ where $\T^1$ is the centralizer of
$\Spin(10)$ in $\Spin(12)$ and the slices becomes $\C \oplus
\C^{16}.$ However, the action on $\C$ is trivial. Then $L$ does
not act coisotropically on $M.$

Since $\C^{27}= \C^{16} \oplus \C^{10} \oplus \C$ as
$\Spin(10)-$submodules, one may prove that $\SU(2) \cdot \T^1
\cdot \Spin (10)$ fails to be multiplicity free. In particular,
see Table {\bf IV}, the subgroups $H$ of $K$ satisfying $\dim H
\geq 47,$ that we have not analyzed yet, are
$$
\SU(2) \cdot \Spin (11) , \ \T^1 \cdot \Spin(11), \ \Spin(11), \
\rho(H)\ H \ \mathrm{simple},\ \rho \in \mathrm{Irr}_{\R}(H),\
\mathrm{d}(\rho)=12.
$$
Let $H= \SU(2) \cdot \Spin (11).$ Since $K \cap \T^1 \cdot \Ea_6 =
\Spin (10)$ then $H \cap \T^1 \cdot \Ea_6 = \T^1 \cdot \Spin
(10),$ so the orbit of $H$ through $[ \T^1 \cdot \Spin (10)  ]$ is
given by $\Spin (11) / \Spin (10) \times \C.$ Note that $H$
preserves the orbit $K[\T^1 \cdot \Ea_6],$ so the slice is given
by $ \R^{10} \oplus \C^{16}, $ on which $\Spin (10)$ acts
diagonally. Let $v \in \R^{10}$ be a unit vector. The orbit is the
unit sphere on $\R^{10}$ and the slice becomes $\R \oplus \C^{16}$
where $\T^1 \cdot \Spin (9)$ acts on $\C^{16}.$ This is the spin
representation, and taking a unit real vector $w,$ the isotropy
group is $\Spin (7)$ and the slice becomes $\R \oplus \R \oplus
\R^{7} \oplus \R^8$ where $\Spin (7)$ acts both on $\R^8$ and on
$\R^7.$ Since $\Spin (7) / {\rm \ G}_2= S^7$ and ${\rm G}_2 /
\SU(3)=S^5,$ the regular isotropy is $\SU (3)$ and the
cohomogeneity is $4.$ So we have $ 4= {\rm rank}(\SU(2) \cdot
\Spin (11))- {\rm rank} (\SU(3)), $ i.e. the action is
multiplicity free. Notice that the slice fails to be polar (see
\cite{Bergmann}). Similarly we may prove that both the $\T^1 \cdot
\Spin (11) -$action and $\Spin(11)-$action  fail to be
multiplicity free Finally, the last case can be excluded by a
straitforward calculation as lemma \ref{zio}.
%
%
\subsection{The case $K=\SU (8)/ \Z_2$}
$K$ is a symmetric group of $\Ea_7$ so $K$ acts coisotropically on
$M.$ We are going to analyze its subgroups. Since $\K \cap \T^1
\cdot \Ea_6 $ is a symmetric group of $K$ and of $\T^1 \cdot
\Ea_6, $ we easily prove that $\K \cap \T^1 \cdot \Ea_6= \T^1
\cdot \SU(2) \cdot \SU(6)$ and the slice becomes $\Lambda^2
(\C^6)$ where $\T^1 \cdot \SU (6) $ acts. Indeed, $K$ is a
symmetric group and the orbit through $[\T^1 \cdot \Ea_6]$ is a
complex orbit so the slice must be a multiplicity-free
representation with degree $15.$ By Tables Ia, Ib and Tables IIa,
IIb we get that the unique possibility is $\Lambda^2 (\C^6).$

By Table {\bf V} and dimensional reasons we may investigate only
$\SUr{1}{7},$ $\SU (7)$ and $\rho(H),$ $H$ is a simple group, such
that $\rho \in Irr_{\C}(H)$ with $d( \rho)=8.$ The last case can
be excluded by a straitforward calculation, while $\SUr{1}{7}$
acts multiplicity free. Indeed, the orbit of $K$ through $[\T^1
\cdot \Ea_6 ]$ is a complex orbit, that is $\SU(8) / \SUr{2}{6},$
the complex Grassmannians of two plane. We may consider the plane
$\pi= \langle e_1, e_2 \rangle$ so the orbit $\SUr{1}{7} \pi$ is
the complex orbit $\SUr{1}{7} / {\rm S}({\rm U}_{1} \times {\rm
U}_1 \times {\rm U}_{6})$ which slice in $M$ is given by $\C^6
\oplus \Lambda^2 (\C^6).$ By Table IIa this action is multiplicity
free. Notice that the slice is not polar. Similarly, one may prove
that also $\SU (7)$ acts coisotropically, but non-polarly, on
$\Ea_7 / \T^1 \cdot \Ea_6 .$
%
%
%
%
%
\section{$M=\Ea_6 / \T^1 \cdot \Spin (10) $}
In this section we analyze the behaviour of the subgroup of
$\Ea_6.$ By dimensional condition, if a subgroup $K\subseteq
\Ea_6$ acts coisotropically on $M=\Ea_6 / \T^1 \cdot \Spin (10) ,$
then $\dim K \geq 26.$ The maximal subgroups of $\Ea_6$ which
satisfy the above inequality (see \cite{Kol} page 41) are the
following
$$
\begin{array}{|l|c|c|c|} \hline
{\rm maximal \ rank}
                     & \T^1 \cdot \Spin(10)
                     & \SU(2) \cdot \Spin (12)
                     & \Sp (1) \cdot \SU (6)   \\ \hline
{\rm no \ maximal \ rank} & \Sp(4) & \Fa_4  &  \\ \hline
\end{array}
$$
%
%
\subsection{The fixed point case $K=\T^1 \cdot \Spin (10) $}
$K$ acts coisotropically and the slice representation appears in
Table Ia and the scalar can be removed. Now, by  Table {\bf IV},
we shall analyze the following cases.
\begin{enumerate}
\item $H=\T^1 \cdot \Spin (k) \times \Spin (10-k).$ Since $\dim H
\geq 26$ we must consider only the cases $\T^1 \cdot \Spin (9),$
$\T^1 \cdot (\T^1 \times \Spin (8))$ and $\T^1 \cdot \Spin(8).$
The first one acts coisotropically but the scalar cannot be
removed. In the other cases, the slice becomes $\C^{16}= \C^8
\oplus \C^8,$ on which $\Spin (8),$ so $\T^1 \cdot (\T^1 \times
\Spin(8))$ acts coisotropically but the scalar cannot be reduced.
Notice that in these cases the slice fail to be polar (see
\cite{Bergmann} and \cite{He}). \item $H=\T^1 \cdot \U(5).$  It is
well know that the isotropy group of $[v]$ in $\P (\C^{16}),$
where $v$ is the highest weight is $\U(5).$ Moreover, the center
of $\U(5)$ acts as scalar while $\SU (5)$ acts trivially on $v.$
Hence $\Spin (10)v= \Spin (10) / \SU (5)$  and the isotropy
representation is given by $\C^5 \oplus \Lambda^2 (\C^5) \oplus
\R.$ In particular $\C^{16}= \C^5 \oplus \Lambda^2 (\C^5) \oplus
\C,$ as $\U(5)-$submodules so by Table IIa this actions is
multiplicity free. Notice that the slice fails to be polar by
Theorem 2 \cite{Bergmann} and  for dimensional reasons any proper
subgroup does not act coisotropically. \item $H= \T^1 \cdot \rho
(H'),$ where $\rho \in Irr_{\C} (H'),$ $d(\rho)=10.$ This case can
be excluded by a straitforward calculation as in lemma \ref{zio}.
\end{enumerate}
%
%
%
\subsection{The case $K=\SU (2) \cdot \SU(6)$}
$K$ acts multiplicity-free since it is a symmetric group of $\Ea_6
.$ We recall that in $\Ea_6$ two involutions $\sigma, \tau$
commuting up to conjugation, i.e. there exists $g \in \Ea_6$ such
that $\sigma$ commutes with ${\rm Ad}(g) \circ \tau \circ {\rm
Ad}(g^{-1})$ (see \cite{co}). In particular we may assume that $K
\cap  \T^1 \cdot \Spin (10)$ is a symmetric group both of $K$ and
of $\T^1 \cdot \Spin (10).$ Hence, looking by the extended Dynkin
diagram of $\Ea_6 ,$ we have $Lie( K \cap  \T^1 \cdot \Spin
(10))=\R + (\R + \su 5 ) \subseteq \asp (1) + \su 6.$ Hence the
orbit through $[\T^1 \cdot \Spin (10)]$ is a complex orbit and the
slice is given by $\Lambda^2 (\C^5).$ Now, we must consider the
maximal subgroup of $K.$ The group $\T^1 \cdot \SU(6)$ acts
coisotropically since the orbit through $[\T^1 \cdot \Spin (10)]$
is $\P (\C^{5})$ and the slice becomes $\C \oplus \Lambda^2 (\C^5
))$ on which $\T^1 \times \U(5)$ acts. In particular $\SU(6)$ does
not act coisotropically since on the slice appears $\C$ on which
the action is trivial. By dimensional condition, one may
investigate  only the following cases: $\T^1 \times \SUr{1}{5}$
and $\T^1 \times \rho (H),$ $H$ simple, $\rho \in Irr_{\C}(H),$
with $d(\rho)=6.$ The second case can be excluded by a
straitforward calculation. In the first case, one may note that
the orbit through
 $[\T^1 \cdot \Spin (10)]$ is a complex orbit and the slice becomes
$\Lambda^2 (\C^5) \oplus \C^5$ where $\U(5)$ acts diagonally.
Hence the slice is a multiplicity free representation which is not
polar by Theorem 2 in \cite{Bergmann}.
%
%
\subsection{The case $K=\Sp(4)$}
$K$ is a symmetric group so the $K-$action is multiplicity free.
By dimensional condition, we shall investigate the cases
$\rho(H),$ $H$ simple, $\rho$ an irreducible representation of
quaternionic type with $d(\rho)=8.$ However, it is easy to check
that this case can be excluded.
%
%
%
\subsection{The case $K={\rm F}_4$}

Since $K$ is a symmetric group the $K-$action is multiplicity
free. Moreover the unique maximal subgroup $H$ which satisfies
$\dim H \geq 26$ is $\Spin (9) \subseteq \Spin (10)$ so we fall
again in the fixed point case.
%
%
%
%
%
\section{Polar actions}
In this section we study which coisotropic actions are polar. It
is well known \cite{PT1} that if a K-action is polar on $M$ then
every slice representation of $K$ is polar. Notice also that the
reducible actions arising from Tables IIa and IIb are not polar;
this can be easily deduced as an application of Theorem $2$ (page
313) \cite{Bergmann}, while see \cite{He} and \cite{Kol}, in the
irreducible case we know  that $\un m$ on $\Sp(m) / \U(m),$ $\un
m$ and $\su m$ when $m$ is odd on $\SO(2m) / \U(m),$ $\Spin (10)$
and $\T^1 \cdot \Spin (10)$ on $\Ea_6 / \T^1 \cdot \Spin (10),$
$\T^1 \cdot \Ea_6$ on $\Ea_7 / \T^1 \cdot \Ea_6$ give rise to
hyperpolar actions. Moreover, any symmetric group and
cohomogeneity one actions are hyperpolar. Hence we may consider
the following cases: $\mf z + \mf t_3$ and $\mf z + \su 2 $ acting
on $\SO(6) / \U(3),$ $\mf z + \asp (2)$ $\asp(1) \otimes \asp(2)$
acting on $\SO(8)/ \U(4),$ $\T^1 \cdot \Spin(12)$ on $\Ea_7 / \T^1
\cdot \Ea_6$ and finally $\asp(m-1) + \un 1$ acting on $\Sp(m) /
\U(m),$ Firstly, we consider $\T^1 \cdot \Spin(12)$ on $\Ea_7 /
\T^1 \cdot \Ea_6.$ We recall that $\T^1$ is not the centralizer of
$\Ea_6$ in $\Ea_7.$ In section $5.1$ we have determined a complex
orbit an its slice is given by $\C \oplus \C^{16}$ on which $\T^1
\cdot (\T^1 \cdot \Spin(10))$ acts. Hence the cohomogeneity is
$3.$ If the action were polar the slice would be a compact
non-flat locally symmetric space. Hence the slice must be a
quotient of $S^3$ and its the tangent space is given by $\R + m,$
where $m$ is a section corresponding to the case $\SU(2) \cdot
\Spin(12),$ so $[m,m]=0,$ since this action is hyperpolar. This
means that the slice has an isometric group of rank at least two,
which is a contradiction.

The case $\asp(1) \otimes \asp(2)$ can be excluded similarly.
Indeed, we have proved that a slice is given by $\C^5$ on which
$\T^1 \cdot \SO(5)$ acts. If the action were polar the section $m$
would be an abelian subspace of dimension $2,$ i.e. the action
would be hyperpolar which is a contradiction, see  \cite{Kol}.

The other cases can be excluded using the same idea. For example,
let $\mf l= \mf z + \su 2 .$ We have proved that the slice
$\Lambda^2 (\C^3 )=\Lambda^2 (\C^2) \oplus (\C \otimes \C^2)^* ,$
so that the action has cohomogeneity $2.$ If the action were polar
a section can be taken as direct sum of the section for the action
of $\T^1$ on $\C$  plus a section for the $\T^1 \cdot \SU (2)$
action on $\Lambda^2 (\C^3) .$ Let $\mf m=<X,Y>,$ where
$$X=\left ( \begin{array}{ccc}
0  & 0    & 0 \\
0  & 0    & 2+i \\
0  & -2-i & 0 \\
\end{array}
\right) \in \Lambda^2 (\C^2), \ \ Y=\left( \begin{array}{ccc}
0   & 1 & 0    \\
-1  & 0 & 0     \\
0   & 0 & 0     \\
\end{array}
\right) \in (\C \otimes \C^2)^* .
$$
One may prove that $[[X,Y],X]$ does not belong to $\mf m.$ Hence,
by Theorem 7.2 page 226 \cite{Hel} on Lie triple system, the
section $\Sigma= \exp ( \mf m )$ is not totally geodesic, hence
the action cannot be polar.
%
%
\section{Appendix}
\newpage
\begin{center}
{{\bf Table I a: } Lie algebras $\mf k$ s.t. $\R+\mf k$ gives rise
to irreducible multiplicity free actions}
\end{center}
$$
\begin{array}{|lrclr|}
\hline {\mf{su}(n)} & n \geq 1 &\quad\quad & {\mf{so}(n)} & n \geq
3
\\
{\mf{sp}(n)} & n \geq 2 & & {S^2(\mf{su}(n)) } & n \geq 2 \\
{\Lambda^2 (\mf{su}(n)) } & n \geq 4 & & {\mf{su}(n) \otimes
\mf{su}(m)} &
n,m \geq 2 \\
{\mf{su}(2) \otimes \mf{sp}(n)} $\ \ \ \ $ & n \geq 2 & &

{\mf{su}(3) \otimes \mf{sp} (n)} $\ \ \ \ $&
n \geq 2 \\
{\mf{su}(n) \otimes \mf{sp}(2)} & n \geq 4
& & {\mf{spin} (7)} &  \\
{\mf{spin} (9)} &  & & {\mf{spin} (10)} &  \\
{\mf{g}_2} & n \geq 1 & & {\mf{e}_6} & n \geq 3 \\ \hline
\end{array}
$$
$\ $ \\
\begin{center}
{{\bf Table I b:} Irreducible coisotropic actions in which the
scalars are removable}
\end{center}
$$
\begin{array}{|lrclr|} \hline
\mf{su} (n) & n \geq 2 & \quad\quad& \asp (n) & n \geq 2 \\
\Lambda^2 ( \mf {su} (n)) & n \geq 4&
 & \mf {su} (n) \otimes \mf {su} (m) & n,m \geq 2,\ n \neq m \\
\mf{spin} (10) & & & \mf {su} (n) \otimes \asp (2) & n \geq 5
\\ \hline
\end{array}
$$
$\ $ \\
\begin{center}
{{\bf Table II a:} Indecomposable coisotropic actions in which the
scalars can be removed or reduced}
\end{center}
$$
\begin{array}{|lr|} \hline
{\mf{su}(n) \mioplus[\mf{su}(n)] \mf{su}(n)}
& n \geq 3,\ a \neq b  \\
{\mf{su}(n)^* \mioplus[\mf{su}(n)] \mf{su}(n)}
& n \geq 3\ a \neq -b \\
{\mf{su}(2m) \mioplus[\mf{su}(2m)] \Lambda^2 (\mf{su}(2m)) }
& m \geq 2, \ b \neq 0  \\
{\mf{su}(2m+1) \mioplus[\mf{su}(2m+1)] \Lambda^2 (\mf{su}(2m+1)) }
& m \geq 2, \ a \neq -mb   \\
{\mf{su}(2m)^* \mioplus[\mf{su}(2m)] \Lambda^2 (\mf{su}(2m)) }
& m \geq 2, \ b \neq 0 \\
{\mf{su}(2m+1)^* \mioplus[\mf{su}(2m+1)] \Lambda^2 (\mf{su}(2m+1))
}
& m \geq 2, \ a \neq mb \\
{\mf{su}(n) \mioplus[\mf{su}(n)] (\mf{su}(n) \otimes \mf{su}(m)) }
& 2 \leq n < m,\ a \neq 0   \\
{\mf{su}(n) \mioplus[\mf{su}(n)] (\mf{su}(n) \otimes \mf{su}(m)) }
& m \geq 2,\  n \geq m+2, a \neq b    \\
{\mf{su}(n)^* \mioplus[\mf{su}(n)] (\mf{su}(n) \otimes \mf{su}(m))
}
& 2 \leq n <m, a \neq 0 \\
{\mf{su}(n)^* \mioplus[\mf{su}(n)] (\mf{su}(n) \otimes \mf{su}(m))
}
& 2 \geq m, n \geq m+2, \ a \neq b \\
{(\mf{su}(2)\otimes \mf{su}(2)) \mioplus[\mf{su}(2)] (\mf{su}(2)
\otimes
\mf{su} (n)) } & n \geq 3, \ a \neq 0 \\
{(\mf{su}(n)\otimes \mf{su}(2)) \mioplus[\mf{su}(2)] (\mf{su}(2)
\otimes \mf{sp} (m)) } & n \geq 3, m\geq 4,  b \neq 0 \\\hline
\end{array}
$$
$\ $ \\
\begin{center}
{{\bf Table II b:} Indecomposable coisotropic actions in which the
scalars cannot be removed or reduced}
\end{center}
$$
\begin{array}{|lr|}  \hline
{\mf{su}(2) \mioplus[\mf{su}(2)] \mf{su}(2)} &\\
{\mf {su}(n)^{(*)} \mioplus[\mf{su} (n)^*] (\mf{su} (n) \oplus
\mf{su} (n)) }
& n \geq 2 \\
{(\mf{su}(n+1)^{(*)} \mioplus[\mf{su}(n+1)] (\mf{su}(n+1) \otimes
\mf{su}(n)) }
& n \geq 2  \\
{(\mf{su}(2) \mioplus[\mf{su}(2)] (\mf{su}(2) \otimes \mf{sp}(m))
}
& m \geq 2 \\
{(\mf{su}(2) \oplus \mf{su}(2)) \mioplus[\mf{su}(2)] (\mf{su}(2)
\otimes \mf{sp}
(m)) } &\\
{(\mf{sp} (n) \oplus \mf{su}(2)) \mioplus[\mf{su}(2)] (\mf{su}(2)
\otimes \mf{sp} (m)) }
& n,m \geq 2 \\
{\mf{sp} (n) \mioplus[\mf{sp} (n)] \mf{sp} (n)}
& n \geq 2 \\
{\mf{spin} (8) \mioplus[\mf{spin} (8)] \mf{so}(8)} & \\ \hline
\end{array}
$$ \\
\begin{small}
\n In the previous Tables we use the notation of \cite{BR}, as an
example $\su n \oplus_{\su n }{\su n }$  denotes the Lie algebra
$\su n $ acting on $\C^n \oplus \C^n$ via the direct sum of two
copies of the natural representation.
\end{small}
$\ $ \\
%
%
\begin{center}
{{\bf Table III:} Maximal subgroups of $\Sp(m)$}
\end{center}
$$
\begin{array}{|r|c|l|} \hline
i)    & \U(m)                    &  \\ \hline ii)   & \Sp(k)
\times \Sp(m-k)    & 1 \leq k \leq m-1  \\ \hline iii)  & \SO(p)
\otimes \Sp(q)    & pq=m,\ p \geq 3,\ q \geq 1  \\ \hline iv)   &
\rho(H)                  & H \ \rm{simple}, \ \rho \in
\Irr_{\H}(H),\ d(\rho)=2m  \\ \hline
\end{array}
$$
%
%
$\ $ \\
\begin{center}
{{\bf Table IV:} Maximal subgroups of $\SO(m)$}
\end{center}
$$
\begin{array}{|r|c|l|} \hline
i) & \SO (k) \times \SO (m-k)& 1 \leq k \leq m-1  \\ \hline ii) &
\SO (p) \otimes \SO (q)& pq=m,\ 3 \leq p \leq q \\ \hline iii) &
\U(k) &  2k=m\\ \hline iv)  & \Sp(p) \otimes \Sp(q) & 4pq=m \\
\hline v)   & \rho(H) & H \ \rm{simple}, \ \rho \in \Irr_{\R}(H)
,\ d (\rho) =m  \\ \hline
\end{array}
$$\\
%
%
\begin{center}
{{\bf Table V:} Maximal subgroups of $\SU(m)$}
\end{center}
$$
\begin{array}{|r|c|l|} \hline
i) & \SO(m) &  \\ \hline ii) & \Sp(n) & 2n=m  \\ \hline iii) &
\SUr{k}{m-k}  & 1 \leq k \leq m-1 \\ \hline iv)  & \SU(p) \otimes
\SU(q) & pq=m,\ p \geq 3,\ q \geq 2  \\ \hline v)   & \rho(H) & H
\ \rm{simple}, \ \rho \in \Irr_{\C}(H),\ d (\rho) =m  \\ \hline
\end{array}
$$
\bibliographystyle{amsplain}

\begin{thebibliography}{10}
\bibitem{Ad}
{\sc J.F. Adams,}
\newblock \textit{Lectures on Exceptional Lie groups},
\newblock  Chigago lectures in Mathematical series, Chicago 1996,
\textbf{181}, 152--186 (1996).
\bibitem{BR}
{\sc C. Benson and G. Ratcliff,}
\newblock \textit{A classification of multiplicity free actions},
\newblock  J. Algebra
\textbf{181}, 152--186 (1996).
%
\bibitem{BA}
{\sc L. Biliotti and A. Gori,}
\newblock \textit{Coisotropic and polar actions on complex Grassmannians,}
\newblock Trans. Am. Math. Soc. \textbf{357}. 1731--1751 (2005).
%
\bibitem{Bergmann}
{\sc I. Bergmann,}
\newblock \textit{Reducible polar representations,}
\newblock  Manus. Math.
\textbf{104}, 309--324 (2001).
%
\bibitem{BTD}
{\sc T. Br\"ocker and T. tom Dieck,}
\newblock \textit{Representation of Compact Lie Groups}
\newblock Springer-Verlag, Berlin-New York, 1985.
%
%
\bibitem{co}
{\sc L. Conlon,}
\newblock \textit{Remarks on commuting involutions,}
\newblock  Proc. Am. Math. Soc.
\textbf{22}, 255--257 (1969).
%
\bibitem{Da}
{\sc J. Dadok J,}
\newblock \textit{Polar cordinates induced by actions of compact Lie Groups,}
\newblock  Trans. Am. Math. Soc.
\textbf{288}, 125--137 (1985).
%
\bibitem{di1}
{\sc E. B. Dynkin,}
\newblock \textit{Semisimple subalgebra of the semisimple Lie algebra,}
\newblock  AMS Translations \textbf {ser. 2}
\textbf{vol 6}, 111--244 (1952).
%
\bibitem{di2}
{\sc E. B. Dynkin,}
\newblock \textit{The maximal subgroups of the classical groups,}
\newblock  AMS Translations, \textbf{ser. 2},
\textbf{vol 6}, 245--278 (1952).
%
\bibitem{He}
{\sc J. Eschenburg and E. Heintze,}
\newblock \textit{On the classification of Polar representation,}
\newblock  Math. Z.
\textbf{286}, 391--398 (1999).
%
%
\bibitem{GS}
{\sc V. Guillemin and S. Sternberg,}
\newblock \textit{Symplectic Techniques in Physics}
\newblock Cambridge: Cambridge University Press, 1984.
%
%
\bibitem{Hel}
{\sc S. Helgason,}
\newblock \textit{Differential Geometry, Lie Groups and Symmetric Spaces}
\newblock  New York-London: Academic Press-inc, 1978.
%
\bibitem{He}
{\sc R. Hermann,}
\newblock \textit{Variational completeness for compact symmetric spaces}
\newblock Proc. Am. Math. Soc. \textbf{11}, 554--546 (1960).
%
\bibitem{HM}
{\sc R. Howe and T. Umeda,}
\newblock \textit{The Capelli identity, the double commutant theorem, and
multiplicity-free actions,}
\newblock  Math. Annalen
\textbf{290}, 565--619 (1991).
%
%
\bibitem{Hw}
{\sc A. T. Huckleberry and T. Wurzbacher,}
\newblock \textit{Multiplicity-free complex manifolds,}
\newblock  Math. Annalen
\textbf{286}, 261--280 (1990).
%
\bibitem{Kac}
{\sc V. G. Kac,}
\newblock \textit{Some remarks on Nilpotent orbit,}
\newblock J. Algebra
\textbf{64}, 190--213 (1980).
%
\bibitem{KN} {\sc S. Kobayashi and K. Nomizu,}
\newblock \textit{Foundations of differential geometry,}
\newblock \textbf{vol. I}, Interscience Publishers, J. Wiley \& Sons, 1963.
%
\bibitem{Kol} {\sc A. Kollross,}
\newblock \textit{A classification of hyperpolar and cohomogeneity one
actions,}
\newblock Trans. Am. Math. Soc.
\textbf{354}, 571--612 (2002).
%
\bibitem{Kr}
{\sc H. Kraft,}
\newblock {\textit \{Geometrische Methoden in der Invariantentheorie,}
\newblock Braunschweig, Wiesbaden: Vieweg 1984.
%
\bibitem{On}
{\sc A. L. Onishchik,}
\newblock \textit{Inclusion relations among transitive compact transformation
groups,}
\newblock Trans. Am. Math. Soc.
\textbf{50}, 5--58 (1966).
%
%
\bibitem{PT1}
{\sc R. S. Palais and C. L. Terng,}
\newblock \textit{A general theory of canonical forms,}
\newblock Trans. Am. Math. Soc.
\textbf{300}, 236--238 (1987).
%
%
\bibitem{PT2}
{\sc  R. S. Palais and C. L.Terng,}
\newblock \textit{Critical point theory and submanifold geometry,}
\newblock Lecture Notes in Mathematics \textbf{1353}, Springer-Verlag,
Berlin-New York, 1988.
%
\bibitem{PoT2}
{\sc F. Podest\`a and G. Thorbergsson,}
\newblock \textit{Polar actions on rank one symmetric spaces,}
\newblock J. Differential Geom.
\textbf{53}, 131-1-75 (2002).
%
\bibitem{PoT}
{\sc F. Podest\`a and G. Thorbergsson,}
\newblock \textit{Polar and {C}oisotropic {A}ctions on {K}\"ahler
{M}anifolds,}
\newblock Trans. Am. Math. Soc.
\textbf{354}, 236--238 (2002).
%
%
\bibitem{SK}
{\sc M. Sato, and T. Kimura,}
\newblock \textit{A classification of irreducible prehomogeneous vector spaces
and their relative invariants,}
\newblock Nagoya Math.  J.
\textbf{65}, 1--155 (1977).
%
\end{thebibliography}

\end{document}